\documentclass[12pt,reqno]{amsart}

\usepackage{amsmath,amssymb,latexsym,amscd,amsthm} 
\usepackage{graphicx}
\usepackage{booktabs}

\newtheorem{theorem}{Theorem}[section]  

\theoremstyle{definition}

\newtheorem{question}[theorem]{Question}

\theoremstyle{remark}

\title{More Colourful Simplices}
\author{Antoine Deza}
\address{Advanced Optimization Laboratory,
Department of Computing and Software,
1280 Main St.~West,
McMaster University,
Hamilton, Ontario,
Canada
L8S 4K1.
}
\email{\{deza,xief\}@mcmaster.ca}
\author[Tamon Stephen]{Tamon Stephen}
\address{Department of Mathematics, Simon Fraser University,
8888 University Drive, Burnaby, British Columbia  V5A 1S6 Canada}
\email{tamon@sfu.ca}
\author{Feng Xie}



\def\bara{B\'ar\'any}

\def\cara{Carath\'eodory}

\def\R{\mathbb{R}}
\def\S{\mathbf{S}}
\def\V{\mathcal{V}}
\def\O{\mathcal{O}}
\def\vecv{\mathbf{v}}
\def\Sph{\mathbb{S}}

\def\csd{\mathbf{depth}}

\def\conv{\operatorname{conv}}

\def\interior{\operatorname{int}}
\def\zero{{\bf 0}}

\begin{document}

\maketitle
\begin{abstract}
We show that any point in the convex hull of
each of $(d+1)$ sets of $(d+1)$ points in general position in $\R^d$
is contained in at least $\lceil {(d+1)^2}/{2} \rceil$
simplices with one vertex from each set.
This improves the known lower bounds for all $d \ge 4$.
\end{abstract}
%
%
\section{Introduction}\label{se:intro}
A point $p \in R^d$ has {\it simplicial depth} $k$
relative to a set $S$ if it is contained in $k$ closed simplices 
generated by $(d+1)$ sets of $S$.  
This was introduced by Liu \cite{Liu90} as a statistical 
measure of how representative $p$ is of $S$, and is
a source of challenging problems in computational geometry -- see for instance \cite{FR05}.
More generally, we consider {\it colourful simplicial depth}, 
where the single set $S$ is replaced
by $(d+1)$ sets, or colours, $\S_1,\ldots,\S_{d+1}$, and the 
{\em colourful} simplices containing $p$ are generated by 
taking one point from each set.

Assuming that the convex hulls of the $\S_i$'s contain
$p$ in their interior, {\bara}'s Colourful {\cara} Theorem
\cite{Bar82} shows that $p$ must be contained in some
colourful simplex.  We are interested in determining 
the minimum number of colourful simplices that can contain
$p$ for sets satisfying these conditions.  That is, we
would like to determine $\mu(d)$, the minimum number of colourful
simplices drawn from $\S_1, \ldots, \S_{d+1}$ that contain $p \in R^d$ 
given that $p \in \interior(\conv(\S_i))$ for each $i$.
Without loss of generality, we assume that the points in 
$\bigcup_i \S_i \cup\{p\}$ are in general position.
Besides intrinsic appeal, 
$\mu(d)$ represents the minimum number of solutions
to the colourful linear programming feasibility problem
proposed in \cite{BO97} and discussed in \cite{DHST06}.

The quantity $\mu(d)$ was investigated in \cite{DHST06}, where it
is shown that $2d \le \mu(d) \le d^2+1$, that $\mu(d)$ is even for odd $d$, and that $\mu(2)=5$.
This paper also conjectures that $\mu(d)= d^2+1$ for all $d \ge 1$.  
Subsequently, \cite{BM06} verified the conjecture for $d=3$
and provided a lower bound of 
$\mu(d) \ge \max(3d, \left\lceil \frac{d(d+1)}{5} \right\rceil)$
for $d \ge 3$,
while \cite{ST06} independently provided a lower bound of
$\mu(d) \ge \left\lfloor \frac{(d+2)^2}{4}\right\rfloor$.
In this note we show:

{\flushleft
{\bf Theorem 1}:
For $d \ge 1$, we have $\mu(d) \ge \lceil \frac{(d+1)^2}{2} \rceil$. 
}

\vspace{2mm}
This strengthens the previously known lower bound for all $d \ge 4$.

\section{Preliminaries}\label{se:setup}
Without loss of generality we can take $p=\zero$.
The sets $\S_1, \ldots, \S_{d+1}$ must each contain at least $(d+1)$
points for $\zero$ to be in the interior of their convex hulls, and
since we are minimizing we can assume they contain no additional
points, i.e.~that $|\S_i|=d+1$ for each $i$.  We assume that
all points are distinct, so no point occurs in two $\S_i$'s,
and $\zero$ is not in any $\S_i$.  We can scale the points 
of the $\S_i$'s so that they lie on the unit sphere $\Sph^{d-1}$:
$\zero$ is in a simplex after scaling if and only if it was
in the simplex before scaling.

We call a set of points drawn from the $\S_i$'s {\it colourful} 
if it contains at most one point from each $\S_i$.  
We call a colourful set of $d$ points which misses $\S_i$
an $\widehat{i}$-{\it transversal}.  Note that $\widehat{i}$-transversals
generate full dimensional pointed colourful cones; we will
say that a transversal {\it spans} a point if the point is
contained in the associated cone.
A key observation is that colourful simplices containing 
$\zero$ are generated whenever the antipode of a point of
colour $i$ is spanned by an $\widehat{i}$-transversal.

In particular, we look at the combinatorial
{\it octahedra}, or cross polytopes, generated by pairs
of disjoint $\widehat{i}$-transversals.
We rely on the topological fact that every octahedron $\Omega$ either
covers all of $\Sph^{d-1}$ with colourful cones, or, every
point $x \in \Sph^{d-1}$ that is covered by colourful cones 
from $\Omega$ is covered by at least two distinct such
cones.  In the case where the points of $\Omega$ form an
octahedron in the geometric sense, these correspond to the
cases where $\zero$ is inside and outside $\Omega$ respectively.
For a proof, see for example the {\em Octahedron Lemma} of \cite{BM06}.
We remark that a given octahedron contains
$2^d$ transversals, though we specify only two disjoint ones
to generate it.

Our strategy for finding distinct colourful simplices is to
begin with a transversal that generates at least one colourful
simplex, and get further points from octahedra that include
this transversal.  We will break into cases based on the number
of colourful simplices generated by the initial transversal
and how many of the octahedra cover $\Sph^{d-1}$.

\section{Proof of the Theorem 1}\label{se:main}
We know that at least one colourful simplex contains $\zero$.
Therefore we have an antipode of colour $(d+1)$ lying in the cone
generated by a $\widehat{d+1}$-transversal $T$.
Without loss of generality we can number the points of 
$\S_1,\S_2,\ldots,\S_d$ so that point $(d+1)$ of $S_i$ is 
included in $T$.  The remaining points of the $S_i$'s can
be numbered arbitrarily.  Let $T_i$ be the set that contains
the points numbered $i$ from $\S_1, \S_2, \ldots \S_{d+1}$.
Then each $T_i$ is a $\widehat{d+1}$-transversal and $T_{d+1}=T$.
Further, the sets $T_1, T_2, \ldots, T_{d+1}$ are pairwise disjoint.
Let $L$ be the set of antipodes of colour $(d+1)$ spanned by $T_{d+1}$,
where $|L|=l>0$.  

\subsection{Points from $d$ octahedra that share a transversal}\label{se:octs}
Now consider the $d$ octahedra $\Omega_1, \Omega_2, \ldots, \Omega_d$
given by pairing $T_i$ with $T_{d+1}$ for $i=1,2,\ldots,d$.
Except for the common transversal $T_{d+1}$, every 
$\widehat{d+1}$-transversal found among the $\Omega_i$'s is distinct.
For each $i$, $\Omega_i$ may or may not cover all of $\Sph^{d-1}$.
Suppose that $b$ of the octahedra cover $\Sph^{d-1}$.
There are $(d+1-l)$ antipodes of colour $(d+1)$ that are not spanned
by $T_{d+1}$, and hence must be spanned by a different transversal
from each of these octahedra.  This gives us a total of
$b(d+1-l)$ distinct simplices containing $\zero$.
Now there remain $(d-b)$ octahedra that do not span all of
$\Sph^{d-1}$.  By the Octahedron Lemma, 
each of the $l$ antipodes spanned by $T_{d+1}$ must also be
spanned by a second transversal from the octahedron generated
by $T_{d+1}$ and $T_i$.
So we find an additional $(d-b)l$ distinct simplices along with
the $l$ simplices generated by the antipodes with $T_{d+1}$ itself.
This brings us to a total of: $l+b(d+1-l)+(d-b)l=(d+1)(b+l)-2bl$ 
distinct colourful simplices containing $\zero$ through this
simple argument.

\subsection{Choice of $T_{d+1}$}\label{se:choice}
In the above argument, $T_{d+1}$ can be any $\widehat{d+1}$-transversal
containing an antipode of colour $(d+1)$.  In the construction of
previous lower bounds, it was noted that if $\csd(\zero)$ is low, 
then there must be a portion of $\Sph^d$ lightly covered by colourful 
cones.  That is to say, if each antipode of colour $(d+1)$ is spanned 
by at least $j$ $\widehat{d+1}$-transversals, then 
$\csd(\zero) \ge j(d+1)$.  We can take $T_{d+1}$ to be a transversal
spanning the least covered antipode.
As we move through the possible values of $i$ in the 
argument of Subsection~\ref{se:octs}, whenever the octahedron
fails to cover $\Sph^{d-1}$ we will see a new cone covering the 
lightly covered antipode.  Hence $(j-1)+b \ge d$.
We thus have that $\csd(\zero)$ is at least
$\max[j(d+1),(d+1)(b+l)-2bl]$ with
$j \ge 1, 1 \le b, l \le d$, and $j+b \ge d+1$.

As long as $l \le \frac{d+1}{2}$, this gives the desired result:
by taking either $j \ge \frac{d+1}{2}$ or $b \ge \frac{d+2}{2}$
we get $\csd(\zero) \ge \frac{d^2+2d+1}{2}$.

\subsection{Single transversals spanning many antipodes}\label{se:many}
This leaves only the case where $\l \ge \frac{d+2}{2}$.
In this situation, we begin with $l$ simplices containing $\zero$
differing only in the $(d+1)$st colour.
We can repeat this exercise for each colour, in which case
we will either find that for each colour $i$,
$l_i \ge \frac{d+2}{2}$, or, for some colour $i$, $l_i \le \frac{d+1}{2}$.
In the latter case, we apply the analysis above to get at least
$\frac{d^2+2d+1}{2}$ distinct simplices containing zero.

If it happens that we get $l_i \ge \frac{d+2}{2}$ for each $i$,
then for each $i$ we have a set $L_i$ of at least $l_i$ antipodes
of colour $i$ which lie on a single $\widehat{i}$-transversal $U_i$.
These generate $(d+1)$ sets $X_1, X_2, \ldots X_{d+1}$ of at least 
$l=\min_i(l_i) \ge \frac{d+2}{2}$ colourful simplices.
There may be some duplication between sets, but we note that
the simplices within each set are distinct and differ only in
the $i$th colour.

We can identify the simplices that make up the $X_i$'s with vectors
in $\{1,2,\ldots,d+1\}^{d+1}$.  We find it helpful to consider them
as vectors in $\R^{d+1}$ unrelated to the initial configuration.

A simplex $\alpha_d$ belonging to a given $X_i$ is represented by a vector
in $\R^{d+1}$ in the following way. The axes correspond to the $d+1$
colours,
and the $q$th coordinate is set to the index in $S_q$ of the point of
colour $q$ of $\alpha_d$.
We recall that the index of points in $S_q$ is set by the arbitrary 
numbering of points of colour $q$ proposed at the beginning
of Section~\ref{se:main}.
 
The vectors associated to the simplices from a given $X_i$ lie on a line
segment in the $i$th coordinate direction.
If a simplex is in both $X_i$ and $X_q$, then the associated vector must
lie at the intersection of the corresponding line segments.

\vspace{1mm}
{\bf Lemma:}
There are at most $d$ {\it duplicate vectors} in the union of the $X_i$'s,
where a vector that is in $k+1$ sets is counted as $k$
duplicate vectors.

{\bf Proof:}
Consider adding the sets iteratively.  
We will say that two sets are in the same {\it component}
if they contain a common point, and extend this to an
equivalence relation.  We remark that each component is 
contained in the topological component formed by taking 
the union of the line segments associated to the $X_i$'s,
but a given topological component will contain multiple
components if the points of intersection of the line segments
are not included in the corresponding $X_i$'s.

We begin with $c=0$ components and $k=0$ duplicate vectors.
Each added set either creates a new component or intersects
$r$ components, producing $r$ duplicate vectors while reducing 
the number of components by $(r-1)$ through the equivalence relation.  
Therefore at each step $c+k$ increases by 1.  
Upon termination, we will have at least
1 component, and hence at most $d$ duplicate vectors.
\qed

\vspace{2mm}
Then the $X_i$'s contain distinct simplices except possibly for up 
to $d+1-c \le d$ repeats arising in this construction, 
where $c$ is the number of components.
This gives us a total of $(d+1)l-(d+1-c)=(d+1)(l-1)+c$ distinct
simplices containing $\zero$.  

However, if $c$ is small, we can readily find additional 
distinct simplices containing $\zero$ by observing that for a
fixed colour $i$, for instance one attaining $l=l_i$, we also
have $(d+1-l)$ antipodes outside of $L_i$.
Each of these antipodes must generate some colourful simplex
containing $\zero$.  
In fact, for each antipode omitted, we could get $\frac{d+1}{2}$
simplices since either $l_i$ or $b$ is this large, but it does not improve
our worst case.
Call this set of simplices $M$, and again consider them as vectors
in $\R^{d+1}$.
They are not included among the vectors associated to simplices 
in $X_i$, since they have different values of coordinate $i$.

The vectors associated to simplices in $M$ could duplicate 
vectors from components other than the one containing $X_i$.  
However, each such component has a fixed value
of colour $i$.  
If $c-1 \ge d+1-l$ it may be the case that all such simplices are
repeats, but our guarantee is $(d+1)(l-1)+c \ge dl+1$.
If $c-1 < d+1-l$ we get at least $d+2-l-c$ additional distinct
simplices from vertices omitted from the $(d+1)$ sets.
This again guarantees us at least
$(d+1)(l-1)+c+(d+2-l-c)=dl+1$ distinct simplices.

\vspace{1mm}
Now as $l \ge \frac{d+2}{2}$ 
we get at least $d\frac{d+2}{2}+1=\frac{d^2+2d+2}{2}$
distinct simplices containing $\zero$.  
Thus our overall worst case for this analysis is 
at $\frac{d^2+2d+1}{2}=\frac{(d+1)^2}{2}$, which can be
rounded up to an integer when $d$ is even.
This improves the known bounds for $d \ge 4$, in particular
from 12 to 13 when $d=4$.
We remark that unlike previous general approaches, this
analysis gives the tight bound of 5 when $d=2$.

\section{A Combinatorial Generalization}\label{se:comb}
The methods in Section~\ref{se:main} rely on the combinatorial
structure of the vectors representing the simplices.
Indeed, there is a nice generalization of the colourful simplicial
depth problem to systems of vectors of in $\{1,2,\ldots,d+1\}^{d+1}$.

Given sets $\S_1, \ldots, \S_{d+1}$ as in Section~\ref{se:intro},
we form the system of vectors $\V$ where $\vecv = (s_1, \ldots, s_{d+1})$
is in $\V$ exactly if the colourful simplex described by $\vecv$
contains $\zero$.  In this context, $\widehat{i}$-{\it transversals} are
simply vectors with the $i$th coordinate removed, and {\it octahedra} are
pairs of disjoint $\widehat{i}$-transversals.  
The system $\V$ has the following two properties:

1. Every element of $\{1,2,\ldots,d+1\}^{d+1}$ is in some $\vecv \in \V$.
This is the combinatorial requirement from {\bara}'s Colourful
{\cara} Theorem.  

2. For any octahedron $\O$, the parity of the set of vectors 
using points from $\O$ and a fixed point $s_i$ for the $i$th 
coordinate is the same for all choices of $s_i$.  
For a system $\V$ arising from colourful simplices, the parity
is odd when the octahedron $\O$ contains $\zero$, and even when it
does not.  This is a purely combinatorial version of the
Octahedron Lemma mentioned in Section~\ref{se:setup}.

\begin{question}\label{qu:gen}
For a given $d \ge 2$, what is the size $\nu(d)$ of a minimal system $\V$ 
of vectors in $\{1,2,\dots,{d+1}\}^{d+1}$ satisfying properties 1 and 2?
\end{question}

The system corresponding to the conjectured minimal core colourful
{\cara} configuration from \cite{DHST06} satisfies properties 1 and 2
with $d^2+1$ vectors, so $\nu(d) \le \mu(d) \le d^2+1$.  
Clearly $\nu(d) \ge d+1$. 
An exhaustive computer search on a laptop shows in a few seconds
that $\nu(2)>4$ and in a few hours that $\nu(3)>8$.
In other words, this approach computationally
verifies that $\mu(2)=5$ and $\mu(3)=10$ (using the fact that $\mu(3)$
must be even).

\section{A Generalized Core}\label{se:gencore}
As a final remark, we mention the recent generalization of the
Colourful {\cara} Theorem in \cite{HPT08} and \cite{AB+09}, in 
which the condition of $\zero$ being in the convex hull of 
each $\S_i$ is relaxed to
require $\zero$ to only be in the convex hull of $\S_i \cup \S_j$
for each $i \ne j$.  
It is natural to ask whether the minimum number of colourful simplices
containing $\zero$ is lower for configurations satisfying 
these weaker conditions.  Call the analogous quantity $\mu^\Diamond(d)$.

In fact, the construction of \cite{DHST06} can be modified in this
to produce configurations showing that $\mu^\Diamond(d) \le d+1$ by
fixing the points of colours $1,2,\ldots,d$ in the same way and then 
clustering all antipodes of the final colour in region that is 
covered by only a single colourful cone from the first $d$ colours.
In this case the relaxed conditions are satisfied almost trivially
since $\zero$ is in $\conv(\S_i)$ for $i=1,2,\ldots,d$.
We note that in this configuration, each colour from $1,\ldots,d$ 
has a unique point which is a generator for all $(d+1)$ colourful
simplices colourful simplices containing $\zero$. 
In other words, in contrast to the situation when $\zero$ is in 
all the $\S_i$'s, some (in fact, most) points from the $\S_i$
generate no colourful simplices containing $\zero$.

The following simple argument shows that $\mu^\Diamond(2)=3$.
Using the assumptions of Section~\ref{se:setup}, we place the
points of the first two colours on the unit circle around $\zero$.
The condition $\zero \in \conv(\S_1 \cup \S_2)$ then means that
every half-circle contains a point from $\S_1 \cup \S_2$.  
If the circle is covered by colourful cones, then each antipode
of the remaining colour generates a colourful simplex containing $\zero$
and we are done.
Otherwise, some segment of the circle is not covered by any colourful cone. 
This segment must be bounded by two points $p$ and $p'$ of 
the same $\S_i$, say $\S_1$.
The three points of $\S_2$ then are on the longer arc between
these points, and for each point of $\S_2$, every point on the 
longer arc is covered by a colourful cone using that point and
either $p$ or $p'$.  The condition that 
$\zero \in \conv(\S_2 \cup \S_3)$ forces at least one of the
antipodes of $\S_3$ to lie in the arc that spans the three
points of $\S_2$.

Finally, we remark that we can generalize $\mu^\Diamond(d)$ 
combinatorially to $\nu^\Diamond(d)$ analogously to Section~\ref{se:comb}.
Combinatorial Property 2 must still hold for such configurations,
but Property 1 fails in the constructions above.  
Nevertheless, we can quickly verify computationally the
$\nu^\Diamond(3)=\mu^\Diamond(3)=4$.

\section{Acknowledgments}
This work was supported by grants from the 
Natural Sciences and Engineering Research Council of Canada (NSERC)
and MITACS, and by the Canada Research Chairs program.
The authors would like to thank the referees for helpful comments
and Imre {\bara} for initiating the discussion in Section~\ref{se:comb}.

\newpage

\newcommand{\etalchar}[1]{$^{#1}$}
\providecommand{\bysame}{\leavevmode\hbox to3em{\hrulefill}\thinspace}
\providecommand{\MR}{\relax\ifhmode\unskip\space\fi MR }
\providecommand{\MRhref}[2]{%
  \href{http://www.ams.org/mathscinet-getitem?mr=#1}{#2}
}
\providecommand{\href}[2]{#2}


\begin{thebibliography}{ABB{\etalchar{+}}09}

\bibitem[ABB{\etalchar{+}}09]{AB+09}
Jorge~L. Arocha, Imre B{\'a}r{\'a}ny, Javier Bracho, Ruy Fabila, and Luis
  Montejano, \emph{Very colorful theorems}, Discrete Comput. Geom. \textbf{42}
  (2009), no.~2, 142--154.

\bibitem[B{\'a}r82]{Bar82}
Imre B{\'a}r{\'a}ny, \emph{A generalization of {C}arath\'eodory's theorem},
  Discrete Math. \textbf{40} (1982), no.~2-3, 141--152.

\bibitem[BM07]{BM06}
Imre B{\'a}r{\'a}ny and Ji{\v r}{\' i} Matou{\v s}ek, \emph{Quadratically many
  colorful simplices}, SIAM Journal on Discrete Mathematics \textbf{21} (2007),
  no.~1, 191--198.

\bibitem[BO97]{BO97}
Imre B{\'a}r{\'a}ny and Shmuel Onn, \emph{Colourful linear programming and its
  relatives}, Math. Oper. Res. \textbf{22} (1997), no.~3, 550--567.

\bibitem[DHST06]{DHST06}
Antoine Deza, Sui Huang, Tamon Stephen, and Tam{\' a}s Terlaky, \emph{Colourful
  simplicial depth}, Discrete Comput. Geom. \textbf{35} (2006), no.~4,
  597--604.

\bibitem[FR05]{FR05}
Komei Fukuda and Vera Rosta, \emph{Data depth and maximal feasible subsystems},
  Graph Theory and Combinatorial Optimization (David Avis, Alain Hertz, and
  Odile Marcotte, eds.), Springer-Verlag, New York, 2005, pp.~37--67.

\bibitem[HPT08]{HPT08}
Andreas~F. Holmsen, J{\'a}nos Pach, and Helge Tverberg, \emph{Points
  surrounding the origin}, Combinatorica \textbf{28} (2008), no.~6, 633--644.

\bibitem[Liu90]{Liu90}
Regina~Y. Liu, \emph{On a notion of data depth based on random simplices}, Ann.
  Statist. \textbf{18} (1990), no.~1, 405--414.

\bibitem[ST08]{ST06}
Tamon Stephen and Hugh Thomas, \emph{A quadratic lower bound for colourful
  simplicial depth}, J. Comb. Opt. \textbf{16} (2008), no.~4, 324--327.

\end{thebibliography}
\end{document}